\documentclass[a4paper,10pt]{article}
\usepackage{amsmath,amsthm}
\usepackage[utf8x]{inputenc}

\newtheorem{propn_hered_limit}{Proposition}
\newtheorem{thm_hered_sat}[propn_hered_limit]{Theorem}
\newtheorem{cor_ps_hered_sat}[propn_hered_limit]{Corollary}
\newtheorem{thm_partial_steiner}[propn_hered_limit]{Proposition}
\newtheorem{lem_edge_partn}[propn_hered_limit]{Lemma}

\title{Supersaturation for hereditary properties}
\author{David Saxton \thanks{Supported by the
  Engineering and Physical Sciences Research Council.}}

\begin{document}



\maketitle

\begin{abstract}
Let $\mathcal{F}$ be a collection of $r$-uniform hypergraphs, and let $0 < p <
1$. It is known that there exists $c = c(p,\mathcal{F})$ such that the
probability of a random $r$-graph in $G(n,p)$ not containing an induced subgraph
from $\mathcal{F}$ is $2^{(-c+o(1)){n \choose r}}$. Let each graph in
$\mathcal{F}$ have at least $t$ vertices. We show that in fact for every
$\epsilon > 0$, there exists $\delta = \delta(\epsilon,p,\mathcal{F}) > 0$ such
that the probability of a random $r$-graph in $G(n,p)$ containing less than
$\delta n^t$ induced subgraphs each lying in $\mathcal{F}$ is at most
$2^{(-c+\epsilon){n \choose r}}$.

This statement is an analogue for hereditary properties of the supersaturation
theorem of Erd\H{o}s and Simonovits. In our applications we answer a question of
Bollob\'as and Nikiforov.
\end{abstract}



\section{Hereditary properties}

Let $\mathcal{F}$ be a collection of $r$-uniform hypergraphs (which we
abbreviate to $r$-graphs). Let $\mathcal{P} = \mbox{Forb}(\mathcal{F})$ be the
collection of all $r$-graphs not containing an induced subgraph from
$\mathcal{F}$. $\mathcal{P}$ is a \emph{hereditary property}; it is a collection
of graphs closed under graph isomorphism and under taking induced subgraphs. Let
$\mathcal{P}^n \subset \mathcal{P}$ be the set of these graphs on $n$ vertices.
Let $G(n,p)$ be a random $r$-graph on $n$ vertices where each edge is included
uniformly and independently with probability $p$.

\begin{propn_hered_limit}[Alekseev~\cite{A1}, Bollob\'as and
Thomason~\cite{BT1}]
 \label{propn_hered_limit}
Let $\mathcal{P}$ be a hereditary property for $r$-graphs. Define $c_n$ via
\[
 \mbox{Pr}[G(n,p) \in \mathcal{P}^n] = 2^{-c_n{n \choose r}}.
\]
Then the limit $\lim_{n \rightarrow \infty} c_n$ exists.
\end{propn_hered_limit}

(Strictly speaking, Alekseev only proved the case $r=2$ and $p=1/2$, but his
argument shows that Proposition~\ref{propn_hered_limit} follows from the
Erd\H{o}s-Hanani conjecture as proved by R\"odl~\cite{RS}. Bollob\'as and
Thomason~\cite{BT1} show that in fact $c_n$ is increasing in $n$.)

For the property $\mathcal{P} = \mbox{Forb}(\mathcal{F})$, let
$c(p,\mathcal{F})$ be the limit in Proposition~\ref{propn_hered_limit}. For
$p=1/2$, the above probability is exactly proportional to the number of graphs
without an induced subgraph from $\mathcal{F}$.

The case $r=2$ has been studied extensively. Pr\"omel and Steger~\cite{PS}
showed that $c(1/2,\{F\}) = 1/t$, where $t$ is the maximum integer for which
there exists $s$, $0 \leq s \leq t$ such that the vertices of $F$ cannot be
partitioned into $s$ cliques and $t-s$ independent sets. For example,
$c(1/2,\{C_4\}) = 1/2$, where $C_4$ is the $4$-cycle.
The case $p \ne 1/2$ has also been studied extensively; for more on this, see
Bollob\'as and Thomason~\cite{BT2} and Marchant and Thomason~\cite{MT}.
The value of $c$ for a collection of $r$-graphs $\mathcal{F}$ when $r>2$ is a
much harder problem; it is at least as hard as the extremal density problem for
$\mathcal{F}$.

\section{Supersaturation}

Our main theorem is that, loosely speaking, the probability of containing a
positive density of induced subgraphs each lying in $\mathcal{F}$ is not that
much less than the probability of containing a single induced subgraph lying in
$\mathcal{F}$.

\begin{thm_hered_sat}
 \label{thm_hered_sat}
Let $\mathcal{F}$ be a collection of $r$-graphs each with at least $t$ vertices,
let $0 < p < 1$ and let $c = c(p,\mathcal{F})$ be as above, i.e., letting
$\mathcal{P}$ be the set of $r$-graphs not containing an induced subgraph lying
in $\mathcal{F}$,
\[
 \mbox{Pr}[G(n,p) \in \mathcal{P}^n] = 2^{(-c+o(1)){n \choose r}}.
\]
Then for every $\epsilon > 0$ there exist $n_0$ and $\delta > 0$ (depending only
on $\epsilon$, $p$, $\mathcal{F}$) such that if $\mathcal{A}$ is a collection of
$r$-graphs on $n > n_0$ vertices with
\[
 \mbox{Pr}[G(n,p) \in \mathcal{A}] > 2^{(-c+\epsilon){n \choose r}},
\]
then some graph in $\mathcal{A}$ contains at least $\delta n^t$ induced
subgraphs each lying in $\mathcal{F}$.
\end{thm_hered_sat}

We wish to draw a parallel with the supersaturation theorem of Erd\H{o}s and
Simonovits~\cite{ES} which says the following. For a collection of forbidden
$r$-graphs $\mathcal{F}$ each on at least $t$ vertices, let
$\mbox{ex}(n,\mathcal{F})$ be the maximum number of edges of an $r$-graph on $n$
vertices containing no copy of any graph in $\mathcal{F}$ (not necessarily
induced), and let $\gamma = \lim_{n \rightarrow \infty} \mbox{ex}(n,\mathcal{F})
{n \choose r}^{-1}$ be the extremal density limit. Then for every $\epsilon > 0$
there exists $\delta > 0$ and $n_0$ (depending only on $\epsilon$ and
$\mathcal{F}$) such that any $r$-graph with more than $(\gamma + \epsilon){n
\choose r}$ edges on $n > n_0$ vertices contains at least $\delta n^t$ subgraphs
each lying in $\mathcal{F}$.

The case $p=1/2$, $r=2$, $\mathcal{F}=\{F\}$ of Theorem~\ref{thm_hered_sat} was
proved by Bollob\'as and Nikiforov~\cite{BN}, using results relying on
Szemer\'edi's Regularity Lemma. They ask whether a proof could be given avoiding
the regularity lemma. We show that this is indeed possible. (More specifically,
for the result they state this can be done by Theorem~\ref{thm_hered_sat}
together with the proof of $c=1/t$ by Alekseev~\cite{A2}.)

Theorem~\ref{thm_hered_sat} also generalizes a theorem of Erd\H{o}s, Rothschild
and Kleitman~\cite{EKR}, where they prove the case $p=1/2$ and
$\mathcal{F}=\{K_t\}$, a complete graph on $t$ vertices.

In~\cite{PS} Pr\"omel and Steger prove that in fact the number of graphs on a
vertex set $V$ of size $|V| = n$ not containing an induced $F$ subgraph is
essentially determined by the number of subgraphs of a single graph. More
specifically, there exist graphs $G = (V,E)$ and $G_0 = (V,E_0)$ with $E \cap
E_0 = \emptyset$ such that every graph $(V,E_0 \cup X)$ with $X \subset E$ does
not contain an induced $F$ subgraph, and the number of generated graphs is
$2^{|E|} = 2^{(1-c){n \choose 2}+o(n^2)}$ with $c = c(1/2,\{F\})$. Write
$\mbox{ex}^*(n,F)$ for the maximum number of edges $|E|$ of such a graph $G$. We
have the following immediate corollary of Theorem~\ref{thm_hered_sat}, which
can be considered a more direct analogue of the supersaturation theorem of
Erd\H{o}s and Simonovits than Theorem~\ref{thm_hered_sat}.

\begin{cor_ps_hered_sat}
 \label{cor_ps_hered_sat}
 Let $F$ be a $2$-graph and let $\epsilon > 0$. Then there exist $n_0$ and
$\delta > 0$ (depending only on $\epsilon,F$) such that for a vertex set $V$ of
size $n = |V| > n_0$, if $G=(V,E)$ is a graph on $\mbox{ex}^*(n,F) + \epsilon {n
\choose 2}$ edges then for every set of edges $E_0 \subset V^{(2)} \setminus E$
there exists a subset $X \subset E$ such that the graph $(V, E_0 \cup X)$
contains at least $\delta n^{|F|}$ induced copies of $F$.
\end{cor_ps_hered_sat}

\section{Proof of Theorem \ref{thm_hered_sat}}

A partial Steiner system with parameters $(r,m,n)$ for a vertex set $V$ of size
$n$ is a collection of sets $\mathcal{D} \subset V^{(m)}$ such that every
$r$-element subset of $V$ appears at most once as a subset of a set in
$\mathcal{D}$. Observe that
\[
 |\mathcal{D}| \leq {n \choose r}  {m \choose r}^{-1}.
\]
As proved by R\"odl~\cite{RS}, there exist partial Steiner systems which cover
almost all $r$-element subsets.
\begin{thm_partial_steiner}[R\"odl]
 \label{thm_partial_steiner} For $r<m$, $\lambda > 0$, there exists $n_0$ such
that for every $n > n_0$, there exists a partial Steiner system $\mathcal{D}$
with parameters $(r,m,n)$ such that
\[
 |\mathcal{D}| \geq (1 - \lambda) {n \choose r} {m \choose r}^{-1}.
\]
\end{thm_partial_steiner}

Let $\mathcal{F},p,\epsilon,c,\mathcal{A},t$ be as in
Theorem~\ref{thm_hered_sat}, and let the common vertex set of the graphs in
$\mathcal{A}$ be $V$.
For a graph $G$ and a subset $D$ of the vertices of $G$, write $G[D]$ for the
induced subgraph on the vertex set $D$. Write also $\mathcal{F} < G$ to denote
that $G$ contains an induced subgraph lying in $\mathcal{F}$.
For a collection of graphs $\mathcal{C}$ on a common vertex set of size $k$,
write
\[
 \mu_k(\mathcal{C}) = \mbox{Pr}[G(k,p) \in \mathcal{C}]
\]
for the measure of the set $\mathcal{C}$ in the space $G(k,p)$.
Thus $ \mu_n(\mathcal{A}) > 2^{(-c+\epsilon){n \choose r}}$.

\begin{lem_edge_partn}
 \label{lem_edge_partn} There exist $\eta, \gamma, \lambda > 0$ and an integer
$m$ (depending only on $\mathcal{F},p,\epsilon$) such that the following is
true. Let $\mathcal{D} = \{D_1, \ldots, D_d\}$ be a partial Steiner system with
parameters $(r,m,n)$ on vertex set $V$ with $d \geq (1 - \lambda) {n \choose r}
{m \choose r}^{-1}$. Let
\[
 I = \{i \in [d] : \mu_n(\{G \in \mathcal{A} : \mathcal{F} < G[D_i]\}) \geq
\gamma \mu_n(\mathcal{A}) \}.
\]
Then $|I| \geq \eta d$.
\end{lem_edge_partn}

\begin{proof}
For $m \geq 1$, let $\mathcal{B}$ be the set of graphs on vertex set $[m]$ that
do not contain an induced subgraph lying in $\mathcal{F}$. Then
$\mu_m(\mathcal{B}) = 2^{(-c+\epsilon'){m \choose r}}$ for some $\epsilon'
\rightarrow 0$ as $m \rightarrow \infty$. Fix $m$ sufficiently large such that
$\epsilon'$ is sufficiently small (to be determined later).

We will choose $\lambda > 0$ later.
Partition $\mathcal{A}$ as $\mathcal{A} = \cup_{S \subset [d]} \mathcal{A}_S$,
where
\[
 \mathcal{A}_S = \{ G \in \mathcal{A} : \{i : \mathcal{F} < G[D_i]\} = S \}.
\]
Let $\theta_i$ be the measure of the set of graphs $G \in \mathcal{A}$ such that
$\mathcal{F} < G[D_i]$, so
\[
 \theta_i = \sum_{S : i \in S} \mu_n(\mathcal{A}_S).
\]
In this notation, $I = \{ i \in [d] : \theta_i \geq \gamma
\mu_n(\mathcal{A})\}$. Let $\eta = |I| / d$. We aim to show that we can take
$\eta > 0$ (independent of $n$). Observe that
\begin{equation}
 \label{ineq_theta}
 \sum_{S \subset [d]} |S| \mu_n(\mathcal{A}_S)=
  \sum_{i \in [d]} \theta_i \leq 
  (\eta d) \mu_n(\mathcal{A}) + ((1 - \eta) \gamma d) \mu_n(\mathcal{A}).
\end{equation}
Observe also that $\mu_n(\mathcal{A}_S) \leq \mu_m(\mathcal{B})^{d-|S|}$ (since
the projection of $\mathcal{A}_S$ onto any $D_i$, $i \not\in S$ is contained
inside a copy of $\mathcal{B}$ on $D_i$). Hence
\begin{align}
 \sum_{S : |S| < \nu d} \mu_n(\mathcal{A}_s) &\leq \sum_{i=0}^{\nu d} {d \choose
i} \mu_m(\mathcal{B})^{d-i} \nonumber \\ 
 &\leq \nu d {d \choose \nu d} \mu_m(\mathcal{B})^{(1-\nu)d} \nonumber \\
 &\leq 2^{O(\nu)d + (-c+\epsilon')(1-\nu)(1-\lambda) {n \choose r}},
\label{eqn_unas_bound}
\end{align}
where $O(\nu) \rightarrow 0$ as $\nu \rightarrow 0$.
Since $\mu_n(\mathcal{A}) = 2^{(-c+\epsilon){n \choose r}}$, we may pick
$\epsilon', \nu, \lambda > 0$ sufficiently small such that the quantity in
(\ref{eqn_unas_bound}) is at most $\mu_n(\mathcal{A})/2$. Thus by
(\ref{ineq_theta}),
\[
 (\eta d) \mu_n(\mathcal{A}) + ((1 - \eta) \gamma d) \mu_n(\mathcal{A}) \geq
  \sum_{S : |S| \geq \nu d} |S| \mu_n(\mathcal{A}_S) \geq
  (\nu d) \mu_n(\mathcal{A}) / 2,
\]
i.e., $\nu / 2 \leq \eta + (1 - \eta) \gamma$. Set $\gamma = \nu / 4$; this
gives $\eta \geq (\nu/4) / (1 - \nu/4) > 0$ as required.
\end{proof}

We are now ready to prove Theorem \ref{thm_hered_sat}.

\begin{proof}
Let $\eta, \gamma, \lambda, m$ be as in Lemma~\ref{lem_edge_partn}. Let $n$ be
sufficiently large for the existence of an $(r,m,n)$ partial Steiner system
$\mathcal{D}$ covering a proportion of at least $1-\lambda$ of the $r$-subsets
of $V$.
Let
\[
 X = \{D \in V^{(m)} : \mu_n(\{G \in \mathcal{A} : \mathcal{F} < G[D]\}) \geq
\gamma \mu_n(\mathcal{A}) \}.
\]

Let $\sigma$ be a randomly and uniformly chosen permutation of $V$, and let
$\mathcal{D}_\sigma$ be the partial Steiner system generated from $\mathcal{D}$
by permuting the vertex set $V$ by $\sigma$. Applying Lemma~\ref{lem_edge_partn}
with $\mathcal{D}_\sigma$ and taking expectations shows that $|X| \geq \eta {n
\choose m}$.

In particular some graph $G \in \mathcal{A}$ contains at least $\eta \gamma {n
\choose m}$ $m$-sets containing an induced subgraph lying in $\mathcal{F}$. Each
fixed copy of an $F \in \mathcal{F}$ is included in at most ${n-t \choose m-t}$
$m$-sets. Hence $G$ contains at least
\[
 \gamma \eta {n \choose m} {n-t \choose m-t}^{-1} \geq \gamma \eta (2m)^{-t} n^t
\]
distinct induced subgraphs each lying in $\mathcal{F}$ (provided $n \geq 2t$).
We can therefore take $\delta = \gamma \eta (2m)^{-t}$, independent of $n$, as
required.
\end{proof}

\paragraph{Acknowledgements} The author would like to thank Andrew Thomason for
his many helpful suggestions.

\end{document}